\documentclass[]{amsart}
\usepackage{amsmath,amsthm,amsfonts, amssymb, color,graphics}
\newtheorem{thm}{\bf Theorem}
\newtheorem{cor}{\bf Corollary}
\newtheorem{lem}{\bf Lemma}
\newtheorem{prp}{\bf Proposition}

\theoremstyle{definition}

\definecolor{wco}{rgb}{0.5,0.2,0.3}

\def\be#1{\begin{equation}\label{#1}}
\def\ee{\end{equation}}

 \def\Re{{\rm Re\,}} \def\Im{{\rm Im\,}\;}\def\res{\mathop{\rm res\,}}
 \def\dist{{\sf dist}}
\def\C{\mathbb{C}}\def\N{\mathbb{N}}\def\Z{\mathbb{Z}}
\def\CL{{\bf C\!\!L}}\def\ds{\displaystyle}\def\ts{\textstyle}
\def\Y{{\mathcal Y}}\def\P{{\mathcal P}}
\def\W{W}
\def\ende{{\sf q.e.d.}}
\def\w{\mathfrak{w}}
\def\z{\mathfrak{z}}\def\y{\mathfrak{y}}

\def\Ende{\quad\ende\medskip}
\def\proof{\noindent{\bf Proof.~}}
\def\remark{\noindent{\bf Remark.~}}
\def\GO#1{O(|z|^{#1})}
\title{Sub-normal Solutions to Painlev\'e's Second Differential Equation}
\author{Norbert Steinmetz}

\begin{document}
\maketitle

\bigskip{\small\begin{quote}{\bf Abstract.} In a recent paper,
Aimo Hinkkanen and Ilpo Laine \cite{HL2} proved that the
transcendental solutions to Painlev\'e's second differential
equation $w''=\alpha+zw+w^3$ have either order of growth
$\varrho=3$ or else $\varrho=\frac32$. We complete this result by
proving that there exist no sub-normal solutions
$(\varrho=\frac32)$ other than the so-called Airy solutions.
\end{quote}


}

\section{\bf Introduction and Main Results}\label{INTRO}

The solutions to the second Painlev\'e differential equation
 $$w''=\alpha+zw+w^3\leqno{[{\rm II}_\alpha]}$$\def\PII{[{\rm
 II}_\alpha]}%
are either rational or transcendental meromorphic functions of
finite order. More precisely, the so-called {\it second Painlev\'e
transcendents} have order of growth $\frac 32\le\varrho\le 3$ (see
Hinkkanen and Laine \cite{HL1}, Shimomura \cite{SSh1, SSh2} and
the author \cite{NSt2,NSt4}). In a recent paper, Hinkkanen and
Laine \cite{HL2} proved that the order is either $\varrho=3$ or
else $\varrho=\frac32$. This result was commonly expected, but
nevertheless marks a great breakthrough.\medskip

The aim of this paper is to describe the solutions of order
$\varrho=\frac32$, called {\it sub-normal}, in more detail, and
using this information to prove the main result on {\it
non-existence} of {\it non-Airy} sub-normal solutions. The
description is intimately associated with the properties of the
{\it first integral}
 $$\W=w^4+zw^2+2\alpha w-w'^2,\quad \W'=w^2.$$
According to \cite{HL2}, the question whether or not $w$ has order
$\varrho=\frac32$ depends on the cluster set $\CL_\varepsilon$ of
the function $\W(z)z^{-2}$ as $z\to\infty$ on
$\C\setminus\P_\varepsilon$. Here $\P$ denotes the set of non-zero
poles of $w$, and $\P_\varepsilon$ denotes the
$\varepsilon$-neighbourhood
 $$\P_\varepsilon={\ts\bigcup_{p\in\P}}\triangle_\varepsilon(p),
 \quad\triangle_\varepsilon(p)=\{z:|z-p|<\varepsilon|p|^{-1/2}\}.$$
The sub-normal solutions are characterised by the conditions
$$n(r,w)=O(r^{3/2})\quad{\rm and}\quad\left\{\begin{array}{rl}{\rm (i)}&\CL_\varepsilon=\{-1/4\}\cr
 {\rm (ii)}&\CL_\varepsilon=\{0\}\end{array}\right. {\rm ~for~ some~} \epsilon>0,$$
and are called of the {\it first} and {\it second kind},
respectively. Special solutions of the first kind are the
so-called Airy solutions, which occur for parameters
$\alpha\in\frac 12+\Z$ and are obtained by (repeated) application
of the so-called B\"acklund transformations to the solutions to
the special Riccati equations
 \be{Airy}w'=\pm(z/2+w^2).\ee

 \begin{thm}\label{HEUREKA}Equation {\rm [II}$_\alpha]$ has no sub-normal solutions other than the
Airy solutions $($which occur for $\alpha\in\frac
12+\Z)$.\end{thm}

The question whether or not zero may be deficient for any
Painlev\'e trans\-cendent is still open. From
$\alpha/w=w''/w-z-2w^2$ and $m(r,w)=O(\log r)$ follows
$m(r,1/w)=O(\log r)$ if $\alpha\ne 0$, hence the value zero is
non-deficient. In case $\alpha=0$ it is well-known and easily
proved that $m(r,1/w)\le \frac12 T(r,w)+O(\log r)$ for any
transcendental solution (see \cite{GLS}, Thm.\ 10.3). As a
by-product of the Hinkkanen-Laine result and Theorem \ref{HEUREKA}
we obtain

\begin{cor}\label{deficiency0}For every solution to $w''=zw+2w^3$
the value zero is non-deficient.\end{cor}

The paper is organised as follows: In section \ref{RESCALE} we
introduce the re-scaling method developed in \cite{NSt2}, which
together with B\"acklund transformations (section \ref{BACKL})
constitutes the main tool. In sections \ref{EINS}, \ref{ZWEI}, and
\ref{RESIDUEN} the solutions of the first and second kind,
respectively, are described in more detail in terms of the
distribution of their poles and residues, while sections
\ref{HEUREKA1}, \ref{HEUREKA2}, and \ref{DEFEKT0} are devoted to
the proofs of Theorem~\ref{HEUREKA} and Corollary
\ref{deficiency0}. Finally, in section \ref{OUTLOOK} we will give
an outlook to sub-normal solutions to Painlev\'e's fourth
equation.

\section{\bf The Re-scaling Method}\label{RESCALE}

The re-scaling method was developed in \cite{NSt2} to prove the
sharp estimate $\varrho\le 5/2$ for the solutions to Painlev\'e's
first equation $w''=z+6w^2$. It also applies to the second and
fourth Painlev\'e equation (see \cite{NSt4}). In the present case,
for any fixed solution to equation {$\PII$} the family
$(w_h)_{|h|\ge 1}$ of re-scaled functions
 $$w_h(\z)=h^{-1/2}w(h+h^{-1/2}\z)$$
is normal in the plane, and every limit function
 \be{LimitFunction}\w=\lim_{h_n\to\infty} w_{h_n}\ee
satisfies
 \be{LimitEq2}\w''=\w+2\w^3,\ee
hence also
 \be{LimitEq}\w'^2=\w^4+\w^2+c.\ee
The constant solutions to (\ref{LimitEq2}) are $\w=0$ and
$\w=\pm\sqrt{-1/2}$, while the non-constant solutions to
(\ref{LimitEq2}) and (\ref{LimitEq}) are either elliptic or
trigonometric functions; the latter only occur in the exceptional cases
$c=1/4$ and $c=0$:
 \be{SPECSOL}\begin{array}{rclr}
 \w&=&\pm\tan(\z/\sqrt2+\tau)/\sqrt2&\quad(c=1/4)\phantom{.}\cr
 \w&=&\pm i/\sin(i\z+\tau)&(c=0).\end{array}\ee
For $c\ne 0, 1/4$ all solutions to equation
(\ref{LimitEq}) occuring as limit functions of the re-scaling
process are elliptic functions.

\medskip{\bf Poles.} The nature of any solution is determined by the
distribution of its poles. The set $\P$ of non-zero poles of some
fixed solution of {$\PII$} is an infinite set, as follows from
$m(r,w)=O(\log r)$ (for notation and results of Nevanlinna Theory
the reader is referred to the monographs of Hayman \cite{WH} and
Nevanlinna \cite{RN}). At any pole $p$ the Laurent series
developments ($\eta=\res\limits_p w=\pm 1$)
 $$\begin{array}{rcl}w(z)&=&\eta (z-p)^{-1}-\frac{1}{6}\eta p(z-p)
 -\frac{1}{4}(\alpha+\eta)(z-p)^2+\mathbf{h}(z-p)^3+\cdots,\cr
 \W(z)&=&-(z-p)^{-1}+10\eta\mathbf{h}-\frac7{36}p^2-\frac13p(z-p)-
 \frac 14(1+\eta\alpha)(z-p)^2+\cdots\end{array}$$
hold; the coefficient $\mathbf{h}$ remains undetermined.
Pre-scribing $\eta$ and $\mathbf{h}$ at $p$ uniquely determines a
solution just like initial values $w_0$ and $w'_0$ at $z_0$ do.
The series converge on some disc $\triangle_\rho(p),$ with
$\rho>0$ independent of $p$.

\medskip{\bf The cluster set}  $\CL_\varepsilon$ is closed and connected (as always), and
also bounded by a constant only depending on $\varepsilon,$ see
\cite{NSt4}, Prop.\ 3.5.
\begin{lem}\label{CLUSTER}The cluster set $\CL=\CL_\varepsilon$
does not depend on $\varepsilon$. Every limit
 \be{POLc}\lim\limits_{p_n\to\infty}[10\eta_n\mathbf{h}_n-\ts\frac7{36}p_n^2]p_n^{-2}
 \quad(\eta_n=\res\limits_{p_n}w),\ee
where $(p_n)$ denotes any appropriate sequence of poles, also
belongs to $\CL.$ Conversely, any limit
$\lim\limits_{h_n\to\infty}h_n^{-2}W(h_n)$ with
 $$\sup\limits_n
|h_n|^{1/2}\dist(h_n,\P)<\infty\quad{\rm and}\quad\inf\limits_n
 |h_n|^{1/2}\dist(h_n,\P)>0$$
coincides with some limit~$(\ref{POLc})$.\end{lem}

\proof The assertions are consequences of the following
observation. If the limit (\ref{LimitFunction}) exists and solves
(\ref{LimitEq}), and if $(k_n)$ denotes any sequence such that
$|h_n|^{1/2}|h_n-k_n|$ is bounded, then some subsequence of
$w_{k_n}$ converges to $\w(\z_0+\z)$ which solves the same
differential equation as does $\w$. \Ende



We note explicitly
$\lim\limits_{p\to\infty}180p^{-2}\mathbf{h}(p)\res_p w=
\left\{\begin{array}{ll}-1 &{\rm if}~ \CL=\{-1/4\}\cr\phantom{-}7&
{\rm if~} \CL=\{0\}\end{array}\right.$, while the solutions to
$w'=\pm(z/2+w^2)$ satisfy $180p^{-2}\mathbf{h}(p)\res\limits_p
w=-1$.

To describe the possible distributions of poles of the second
Painlev\'e transcendents of order $\varrho=\frac32$ we need the
following result on the local distribution of poles; it is based
on the distribution of poles of the limit functions
$\w=\lim\limits_{p_n\to\infty}w_{p_n}$ ($p_n\in\P$).

\begin{lem}\label{LOCDIS}Suppose that $w$ solves $\PII$ and has order of growth $\varrho=\frac32$. Then given
$\varepsilon>0$ and $R>0$ there exists $r_0>0$ such that for any
pole $p$ satisfying $|p|>r_0$, the poles of $w$ in
$\triangle_R(p)=\{z:|z-p|<R|p|^{-1/2}\}$ may be labelled in such a
way that $p_0=p$ and
 $$|p_k-(p+k\sqrt{2}\pi p^{-1/2})|<\varepsilon|p|^{-1/2}
 \quad(-k_1\le k\le k_2)\leqno{{\rm(first~kind)}}$$
and
 $$|p_k-(p+k\pi ip^{-1/2})|<\varepsilon|p|^{-1/2}
 \quad(-k_1\le k\le k_2),\leqno{{\rm(second~kind)}}$$
respectively, hold.\end{lem}

The {\bf proof} is an immediate implication of the re-scaling
method and the known distribution of poles of the solutions
(\ref{SPECSOL}) to the special re-scaled differential equation.

\medskip To determine the asymptotics of the solutions of order
$\varrho=\frac32$ more precisely, we shall repeatedly apply the
following estimates; the first one is an immediate corollary of
the Cauchy integral formula.

\begin{lem}\label{cauchy}Suppose that $f$ is holomorphic in some sector $S:a<\arg z<b$
satisfying $f(z)=\GO{\lambda}$ as $z\to\infty$ in $S$. Then
$f^{(k)}(z)=\GO{\lambda-k}$ as $z\to\infty$ holds in every smaller
sector $S(\delta):a+\delta<\arg z<b-\delta$.\end{lem}

\begin{lem}\label{mittagleffler}Let $(c_k)$ be any complex sequence
$(0<|c_1|\le |c_2|\le\cdots\le |c_k|\to\infty)$ with counting
function $n(r)= {\rm card}~\{c_k:|c_k|\le r\}=O(r^\varrho)$
$(\varrho=h+\gamma,$ $h\in \N_0,$ and $0<\gamma<1).$ Then
$|z-c_k|\ge\kappa\max\{|z|,|c_k|\}$ for some $\kappa>0$ and every
$k$ implies
 $$\sum_{k=1}^\infty\Big|\frac{z^h}{(z-c_k)c_k^h}\Big|=O(|z|^{\varrho-1})\quad(z\to\infty).$$
 \end{lem}

\proof From $n(r)=O(r^\varrho)$ and
 $\ds\Big|\frac{z^h}{(z-c_k)c_k^h}\Big|\le\kappa^{-1}
 \frac{r^{h-1}}{|c_k|^h}\min\Big\{1,\frac{r}{|c_k|}\Big\}$ on $|z|=r$
follows
 $$\sum_{k=1}^\infty\Big|\frac{z^h}{(z-c_k)c_k^h}\Big|\le\kappa^{-1}r^{h-1}
 \int_{|c_1|}^{r}\frac{dn(t)}{t^h}+
 \kappa^{-1}r^{h}\int_{r}^{\infty}\frac{dn(t)}{t^{h+1}}=O(r^{\varrho-1}).\quad\ende$$

\section{\bf Solutions in the Yosida Class}\label{YOSIDACLASS}

If the cluster set $\CL$ contains none of the values $0,-1/4$,
then all limit functions (\ref{LimitFunction}) are non-constant,
hence $w$ belongs to the {\it Yosida Class} $\Y_{\frac 12,\frac
12}$, being defined and discussed in \cite{NSt3}. These solutions
are traditionally called {\it non-truncated} (see Boutroux
\cite{PB1,PB2}). Among others it follows that $T(r,w)\asymp r^3$
and that the poles are {\it regularly distributed}:  given $R>0$
there exists $r_0>0$ and $C>1$, such that any disc
$\triangle_R(z_0)$ with $|z_0|>r_0$ contains at least $C^{-1}R^2$
and at most $CR^2$ poles. This holds in a modified form if the
cluster set is restricted to some sector $S=\{z:\theta_1\le\arg
z\le \theta_2\}$: the poles in $S$ are regularly distributed, and
again $T(r,w)\asymp r^3$ holds.

\section{\bf Solutions of the First Kind}\label{EINS}

Throughout this section $w$ will denote a transcendent of the
first kind.

\medskip{\bf Strings of poles of the first kind.}
A {\it string} in the truncated sector
 $$S'_0: |\arg z|\le\pi/3,~\Re z\ge c_1>0$$
is a sequence $(p_k)_{k=0,1,2,\ldots}$ such that (we assume $\Re
p_k^{-1/2}>0$)
 $$p_{k+1}=p_k+\sqrt{2}\pi p_k^{-1/2}(1+o(1))\quad(k\to\infty);$$
$p_0$ is called the {\it root} of the string
$(p_k)_{k=0,1,2,\ldots}$.

\begin{prp}\label{CHAINEX}For $c_1>0$ sufficiently large, every pole $p_0$ in $S'_0$
is the root of some uniquely determined string of poles
$(p_k)_{k=0,1,2,\ldots}$ contained in $S'_0$. It has the following
properties:
\begin{itemize}\item[a.] $w$ has constant residues on the string;
\item[b.] $\lim\limits_{k\to\infty}\arg p_k=0$;\item[c.] The
counting function of the string satisfies $n(r)=\frac{\sqrt
2}{3\pi}r^{3/2}(1+o(1)).$\end{itemize}\end{prp}

\proof The construction of the sequence $(p_k)$ is obvious. We
denote by $c_n$ the number $r_0$ in Lemma~\ref{LOCDIS} which
belongs to $\varepsilon=\varepsilon_n=10^{-n}$ and
$R=5>\sqrt{2}\pi$, and start with $p_0\in S'_0$. If $p_k$ is
already constructed, then $p_{k+1}$ is uniquely determined by
Lemma~\ref{LOCDIS}. We have, however, to ensure that the procedure
does not break down, that is we have to show that $p_{k+1}\in
S'_0$. Writing $p_k=|p_k|e^{i\theta_k}$ it follows that
 $$\begin{array}{rcl}
 \Re p_{k+1}&>&\Re p_k+|p_k|^{-1/2}(\sqrt{2}\pi\cos(\theta_k/2)-\varepsilon_1)\cr
 &\ge& \Re p_k+(\sqrt{6}\pi/2-\varepsilon_1)|p_k|^{-1/2}>\Re p_k+3|p_k|^{-1/2}.\end{array}$$
Similarly,
 $$|\Im p_{k+1}|\le |\Im
 p_k|+|p_k|^{-1/2}(\varepsilon_1-\sqrt{2}\pi|\sin(\theta_k/2)|)<|\Im p_k|+|p_k|^{-1/2}\varepsilon_1$$
holds. With the help of $\ds\frac{a+b}{c+d}\le\max\Big\{\frac ac,
\frac bd\Big\}$ for $a,b,c,d>0$ we obtain
 \be{winkel}|\theta_{k+1}|=\arctan\frac{|\Im p_{k+1}|}{\Re p_{k+1}}
 \le\arctan\max\Big\{\frac{|\Im p_{k}|}{\Re p_{k}},\frac{\epsilon_1}3\Big\}
 \le \max\{|\theta_k|, \epsilon_1\}.\ee
It is obvious that $\Re p_k\to\infty$ monotonically, and that the
sequence $(|\Im p_k|)$ decreases as long as $|\theta_k|\ge
\varepsilon_1$. From $\ds |\theta_k|<\frac{|\Im p_k|}{\Re p_k}$,
however, follows that $|\theta_{k_1'}|<\varepsilon_1$ for some
$k_1'$, hence $|\theta_k|<\varepsilon_1$ for $k\ge k'_1$ follows
from (\ref{winkel}). If we denote by $k_n$ the first index such
$\Re p_k>c_{n}$, then the above argument shows that there exists
some $k_n'\ge k_n$, such that $|\theta_k|<\varepsilon_n$ holds for
$k>k_n'$. This yields b. To prove c. we consider the conjugate
sequence $q_k=p_k^{3/2}$. From $p_{k+1}=p_k+\sqrt{2}\pi
p_k^{-1/2}+o(|p_k|^{-1/2})$ follows $q_{k+1}=q_k+\frac
32\sqrt{2}\pi+o(1)$, thus $q_k=\frac 32\sqrt{2}\pi k(1+o(1))$,
 $$\ts p_k=\big(\frac
32\sqrt{2}\pi\big)^{2/3}\, k^{2/3}(1+o(1))\quad{\rm and}\quad
n(r)=\frac{\sqrt 2}{3\pi}\,r^{3/2}(1+o(1)).\quad\ende$$

\remark For $p_{-1}\in  S'_0$ the string just constructed may be
uniquely extended ``to the left'' such that $(p_k)_{k>
-k_0}\subset S'_0$ and $p_{-k_0}\notin S'_0$. Relabelling this
string we may thus always assume that
$(p_k)_{k=0,1,2,\ldots}\subset S'_0$, but $p_{-1}\notin S'_0$.
Such a string is called {\it maximal.}

\medskip
There is just one step from local to global distribution of poles.

\begin{thm}\label{GLOBDIS1}Let $w$ be any second Painlev\'e transcendent of the first kind.
Then up to finitely many the poles of $w$ form a finite number
$\ell(w)$ of maximal strings $\sigma=(p_k)_{k=0,1,2,\ldots}$ with
total counting function
 $$n(r,w)=\ell(w)\,\frac{\sqrt 2}{3\pi}\,r^{3/2}(1+o(1)),$$
and such that the following is true:
\begin{itemize}\item $w$ has constant residues on $\sigma$;
\item $\sigma$  is asymptotic to one of the rays $\arg z=0,$ $\arg
z=\frac 23\pi,$ and $\arg z=-\frac 23\pi;$ \item $\sigma$ is
accompanied by a string $(q_k)$ of zeros
$q_k=p_k+\frac{\sqrt2}2\pi p_k^{-1/2}(1+o(1));$ \item any two
strings $(p_k)$ and $(p'_k)$ are separated from each other, i.e.,
$$\lim\limits_{k\to\infty}|p_k|^{1/2}\dist(p_k,\{p'_n\})=\infty.$$
\end{itemize}
Furthermore, $w$ has Nevanlinna characteristic
 $$T(r,w)=\ell(w)\,\frac{\sqrt 8}{9\pi}\,r^{3/2}(1+o(1)),$$
and satisfies $w(z)\sim \sqrt{-z/2}$ as $z\to\infty$ on every
sector $S_0''(\delta):|\arg z-\pi|<{\frac\pi3-\delta}$ and $S_{\pm
1}''(\delta)=e^{\pm 2\pi i/3}S_0''(\delta),$ for some suitably
chosen branch of the square-root depending on the sector.\end{thm}

\proof From $n(r,w)=O(r^{3/2})$ follows that there are only
finitely many strings of poles. This yields $T(r,w)=N(r,w)+O(\log
r)=\ell(w)\frac{\sqrt 8}{9\pi}\,r^{3/2}(1+o(1))$. The asymptotics
for $w$ follows from the fact that the re-scaling process for any
sequence $(h_n)$ with $|h_n|^{1/2}\dist(h_n,\P)\to\infty$ leads to
the limit functions $\w=\sqrt{-1/2}$. \Ende

\medskip{\bf Series expansion.} In \cite{NSt4} it was shown that
for every second transcendent with $w(0)\ne\infty$
 $$w(z)=w(0)+\lim_{r\to\infty}\sum_{|p|\le
 r}\frac{\eta(p)z}{(z-p)p}\quad(\eta(p)=\res\limits_{p} w)$$
holds; if $w$ has a pole at $z=0$, the term $w(0)$ has to be
replaced by $\eta(0)/z$. In our case the above Mittag-Leffler
expansion exists not only as a Cauchy principal value, but
converges absolutely. Then also
 \be{WSERIES}\W(z)= Q(z)-\frac{|\eta(0)|}z-\sum_{p\in\P}\frac{z}{(z-p)p}\ee
holds, where $Q$ is a polynomial of degree at most two (see
\cite{NSt4}, Thm. 4.3). Lemma~\ref{mittagleffler} applies to
$W-Q$, and from $|\W(z)-Q(z)|=O(|z|^{1/2})$ as $z\to\infty$ in
each sector $S_j''(\delta)$ and $\CL=\{-1/4\}$ then follows
$Q(z)=-\ts\frac 14z^2+a_1z+a_0.$ Also in each sector
$S_j''(\delta)$ we get
 $$\begin{array}{rcl}\W&=&-\frac 14z^2+a_1z+\GO{1/2}\cr
 zw^2=zW'&=&-\frac 12z^2+a_1z+\GO{1/2}\cr
 w^4&=&\phantom{-}\frac 14 z^2-a_1z+\GO{1/2}\end{array}$$
 hence $zw^2+w^4-\W=-a_1z+\GO{1/2}$  and
 $w'^2-2\alpha w=\GO{1/2}$ yield $a_1=0$. We have thus proved

 \begin{thm}\label{asymp1}Any first kind transcendent $w$
 satisfies
 \be{wasymp1}w=\sqrt{-z/2}+\GO{-1}\ee
as $z\to\infty$ in every sector $S_j''(\delta)$ $($for some branch
of $\sqrt{-z/2}$, depending on the
 sector$)$,
 $$\W=-{\ts\frac 14}z^2+\GO{1/2}\quad{\rm and}\quad
Q(z)=\ts -\frac 14z^2+a_0.$$
\end{thm}

\section{\bf Solutions of the Second Kind}\label{ZWEI}

Now $w$ will denote a sub-normal solution of the second kind.
Again (\ref{WSERIES}) holds, where now $\deg Q\le 1$ follows from
$\CL=\{0\}$ and Lemma \ref{mittagleffler}. Since by \cite{NSt4},
Thm.\ 4.5, the order of $w$ is $\varrho\ge 2$ (hence $\varrho=3$)
if $\deg Q=1$, we have $\deg Q=0$.

\medskip A {\it string} of poles $(p_k)_{k=0,1,2,\ldots}$ in the sector
 $$S''_0:|\arg z-\pi|<\pi/3,\quad\Re z<-c_1,$$
is now characterised by the condition ($\Im p_k^{-1}\ge 0$)
 $$p_{k+1}=p_k+i\pi p_{k}^{-1/2}(1+o(1))\quad(k=0,1,2,\ldots);$$
it is called {\it maximal} in $S''_0$ if $p_{-1}\notin S''_0$.
Similarly we define strings of poles in the sectors $S''_{\pm
1}=e^{\pm 2\pi i/3}S''_0$, and obtain the following analog to
Theorem \ref{GLOBDIS1}.

\begin{prp}\label{GLOBDIS2}Let $w$ be any sub-normal solution
of the second kind. Then up to finitely many the poles of $w$ form
a finite number $\ell(w)$ of maximal strings
$\sigma=(p_k)_{k=0,1,2,\ldots}$ with total counting function
 $n(r,w)=\ell(w)\,\frac{2}{3\pi}\,r^{3/2}(1+o(1)),$
and such that the following is true:
\begin{itemize}\item the
residues {\it alternate}, i.e.,
$\ds\res\limits_{p_{k+1}}w=-\res\limits_{p_k} w$; \item $\sigma$
is asymptotic to one of the rays $\arg z=\pi$, $\arg z=\frac
13\pi,$ and $\arg z=-\frac 13\pi;$ \item any two strings $(p_k)$
and $(p'_k)$ are separated from each other.
\end{itemize}
 \end{prp}






\section{\bf B\"acklund Transformations}\label{BACKL}

The so-called {\it Airy solutions} are obtained from the solutions
to any of the special Riccati equations (\ref{Airy}) by successive
application of so-called B\"acklund transformations.

Generally spoken, a {\it B\"acklund transformation} is a change of
variables $w_1(\zeta)=w(z)$, $\zeta=az,$ that transforms equation
{$\PII$} into itself or into some equation [II$_{\alpha_1}$]  with
different parameter. Simple examples are $w_1(z)=-w(z)$
($\alpha_1=-\alpha$) and $w_1(z)=\mu w(\mu z)$ ($\mu^3=1$,
$\alpha_1=\alpha$). More sophisticated B\"acklund transformations
are
 \be{SpecBack}w_1=-w-\frac{\alpha+1/2}{w'+w^2+z/2}\quad{\rm and}
 \quad w_{-1}=-w+\frac{\alpha-1/2}{w'-w^2-z/2},\ee
which change $\alpha$ to $\alpha_1=\alpha+1$ and
$\alpha_1=\alpha-1$, respectively. It is obvious that B\"acklund
transformations (\ref{SpecBack}) preserve the order $\varrho$,
and, by Theorem \ref{GLOBDIS1} and Proposition \ref{GLOBDIS2},
even preserve the first and second kind solutions.

\medskip{\bf A special B\"acklund transformation.} In \cite{GLS}
the authors describe the connection between equations [II$_0$] and
[II$_{\frac12}$]. If $y$ is a non-trivial solution to [II$_0$],
then
 \be{trans0}w(z)=-\frac d{dz}\log y(-2^{-1/3}z)\ee
solves [II$_{\frac 12}$] and is not an Airy solution ($w'\ne
z/2+w^2$); conversely, if $w$ solves [II$_{\frac12}$] and is not
an Airy solution, then the function $y$, being defined locally by
 \be{trans1}y^2(-2^{-1/3}z)=-2^{1/3}(w'(z)-z/2-w^2(z))\ee
is a non-trivial solution to [II$_0$]. The poles and zeros of $y$
correspond to poles of $w$ with residues $1$ and $-1$,
respectively. The Airy solutions to [II$_{\frac 12}$] correspond
to the trivial solution $y=0$. It is obvious that this
transformation preserves the order of growth, but interchanges the
transcendents of the first kind (minus the Airy solutions) and
those of the second kind. This follows at once from the
asymptotics of the involved functions, and also from the
distribution of their poles: noting that $\tilde p=ap$ implies
$\tilde p+b\tilde p^{-1/2}=a(p+a^{-3/2}bp^{-1/2})$, we obtain 
 $$a^{-3/2}b=\left\{\begin{array}{rll}\pi i&{\rm if}~ a=-2^{1/3} &{\rm ~and~} b=\sqrt2\pi\cr
\sqrt{2}\pi &{\rm if}~ a=-2^{-1/3} &{\rm ~and~} b=\pi i.
 \end{array}\right.$$
The possible distribution of poles $\oplus$ and $\ominus$ with
residues $1$ and $-1$, respectively, and zeros $\circledcirc$ of a
second kind solution to $y''=zy+2y^3$ (left), and the distribution
of poles and zeros of the corresponding first kind solution
$w=-\frac d{dz}\log y(-2^{-1/3}z)$ (right) along the real axis is
displayed below. The strings are separated from each other.

 $$\begin{array}{c}\cdots\oplus\quad\ominus\quad\oplus\quad\ominus\quad\oplus\quad\ominus
\qquad\qquad\oplus\circledcirc\oplus\circledcirc\oplus\circledcirc\oplus\circledcirc\oplus\circledcirc\oplus\circledcirc\cdots\cr\cr
\cdots\circledcirc\quad\circledcirc\quad\circledcirc\quad\circledcirc\quad\circledcirc\quad\circledcirc
\qquad\qquad\ominus\circledcirc\ominus\circledcirc\ominus\circledcirc
\ominus\circledcirc\ominus\circledcirc\ominus\circledcirc\cdots\cr\cr\end{array}$$

\section{\bf The Distribution of Residues}\label{RESIDUEN}

We shall denote by $n_\oplus(r)$ and $n_\ominus(r)$ the counting
function of poles with residues $1$ and $-1$, respectively. For
second kind transcendents the residues are equally distributed in
each string of poles, hence $n_\oplus(r)-n_\ominus(r)=o(r^{3/2})$
holds.

\medskip Let $w$ be any first kind transcendent. If the circle $|z|=r$ contains no poles, then
 $$\frac{1}{2\pi
i}\int_{|z|=r}w(z)\,dz=n_\oplus(r)-n_\ominus(r)=(\ell_\oplus
-\ell_\ominus)\frac{\sqrt
 2}{3\pi}r^{3/2}(1+o(1))$$
holds, where $\ell_\oplus$ and $\ell_\ominus$ count the number of
maximal strings with residues $+1$ and $-1$, respectively. We
choose $\delta>0$ sufficiently small and replace any arc of
$|z|=r$ that intersects some disc $\triangle_\delta(p)$ by a
sub-arcs of $\partial\triangle_\delta(p)$ (such that
$|w(z)|=\GO{1/2}$) to obtain a simple closed curve $\Gamma_r$.
Then also
 \be{anzahl}\frac{1}{2\pi
 i}\int_{\Gamma_r}w(z)\,dz=n_\oplus(r)-n_\ominus(r)\ee
holds. If $\gamma_r$ and $\gamma_r'$ denote the part of $\Gamma_r$
in $0\le\arg z<\delta$ and $\delta\le \arg z\le\frac 23\pi$,
respectively, then
 $$\Big|\frac{1}{2\pi i}\int_{\gamma_r}w(z)\,dz\Big|<K\delta r^{3/2}$$
 and
 $$\frac{1}{2\pi i}\int_{\gamma'_r}w(z)\,dz =
 \frac{1}{2\pi i}\int_{\gamma'_r}\sqrt{-z/2}\,dz+o(r^{3/2})=\mu\frac{\sqrt2}{3\pi}r^{3/2}+O(\delta r^{3/2}),$$
hold, with $\mu=\pm 1$ depending on the branch of $\sqrt{-z/2}$.
This yields

\begin{prp}\label{LPLUSMINUS}
Any sub-normal solution of the first kind satisfies
$$|\ell_\oplus-\ell_\ominus|=1\quad {\rm or~ else}\quad|\ell_\oplus-\ell_\ominus|= 3.$$\end{prp}

\remark The following results deduced from Proposition
\ref{LPLUSMINUS} for the solutions to the Riccati equation
 \be{Airyplus}w'=z/2+w^2\ee
are well known, see, e.g. \cite{GS}. Since all residues equal $-1$
we have $\l_\oplus(w)=0$, hence either $\ell_\ominus(w)=1$ or else
$\ell_\ominus(w)=3$. There exist three distinguished solutions
$w_1$, $w_{2}(z)=e^{2\pi i/3}w_1(ze^{2\pi i/3})$, and
$w_{3}(z)=e^{-2\pi i/3}w_1(ze^{-2\pi i/3})$ with $\ell_\ominus=1$.
The labelling is chosen in such a way that
$w_1(z)\sim\psi(z)=\sqrt{-z/2}$ with $\Im\psi(z)>0$ holds on
$0<\arg z<2\pi.$ By symmetry and uniqueness, the poles of $w_1$
are real and positive. For any solution $w_0\ne w_k$ to
(\ref{Airyplus}) we have
 $$w_0(z)\sim\left\{\begin{array}{lr}\phantom{-}\psi(z)&\qquad(0<\arg<\frac
23\pi)\phantom{.}\cr -\psi(z)&(\frac 23\pi<\arg<\frac
 43\pi)\phantom{.}\cr \phantom{-}\psi(z)&(\frac 43\pi<\arg<2\pi).\end{array}\right.$$

 \section{\bf Proof of Theorem \ref{HEUREKA}: First Kind Solutions}\label{HEUREKA1}

Let $w$ be any sub-normal solution of the first kind to
[II$_\alpha$]. We first assume $\alpha\notin\frac 12+\Z$, and set
$V=w'+w^2+z/2$, $w_1=-w-\ds\frac{\alpha+1/2}V$, and
 $$\Delta(w)=\ell_\oplus-\ell_\ominus.$$
Then $w_1$ solves [II$_{\alpha+1}$], and the poles of $w$ and
$w_1$ and the zeros of $V$ are related as follows:

\be{POLE}
\begin{array}{rl}{\rm(i)}& \res\limits_p w=-1\Rightarrow V(p)=\infty {\rm ~(doubly)~and~}\res\limits_p
w_1=1;\cr {\rm(ii)}& \res\limits_p w=1\Rightarrow V(p)=0 {\rm
~and~} \res\limits_pw_1=0;\cr {\rm(ii)} &\res\limits_p w=0 {\rm
~and~} V(p)=0\Rightarrow\res\limits_p w_1=-1.
\end{array}\ee

The distribution of zeros $\circledcirc$ of $V$ (left and right),
and poles $\oplus$ and $\ominus$ with residues $1$ and $-1$,
respectively, of both $w$ (left) and $w_1$ (right):

\be{VERTEILUNG}\begin{array}{rccc}{\rm(i)}&\ominus\quad\ominus\quad\ominus\quad\ominus\quad\ominus\quad\ominus&\qquad&
\oplus\quad\oplus\quad\oplus\quad\oplus\quad\oplus\quad\oplus\cr
&&&\cr
{\rm(ii)}&\oplus\quad\oplus\quad\oplus\quad\oplus\quad\oplus\quad\oplus&&\circledcirc\quad\circledcirc\quad
\circledcirc\quad\circledcirc\quad\circledcirc\quad\circledcirc\cr
&&&\cr {\rm(iii)}&\circledcirc\quad\circledcirc\quad
\circledcirc\quad\circledcirc\quad\circledcirc\quad\circledcirc&\qquad&
\ominus\quad\ominus\quad\ominus\quad\ominus\quad\ominus\quad\ominus\end{array}\ee

First of all we obtain $\ell_\oplus(w_1)=\ell_\ominus(w)$ from
(i), while
 $$m(r,1/V)\le m(r,w_1)+m(r,w)+O(1)=O(\log r),$$
hence
 $$\begin{array}{rcl}N(r,1/V)&=&N(r,V)+O(\log r)=2N_\ominus(r,w)+O(\log r)\cr
 &=&N_\ominus(r,w)-[N_\oplus(r,w)-N_\ominus(r,w)]+N_\oplus(r,w)+O(\log r)\end{array}$$
and (ii) and (iii) imply
 $$\ell_\ominus(w_1)=\ell_\ominus(w)-\Delta(w)\quad{\rm and}\quad\Delta(w_1)=\Delta(w).$$
Repeated application \big($w_{\nu+1}=-w_{\nu}-\frac{\alpha+\nu+
1/2}{w_\nu'+w_\nu^2+z/2}$\big) yields
 $$\begin{array}{rcl}\ell_\ominus(w_\nu)&=&\ell_\ominus(w)-\nu\Delta(w),
 \quad\Delta(w_\nu)=\Delta(w),\quad{\rm and}\cr
 \ell_\oplus(w_\nu)&=&\ell_\ominus(w_{\nu-1})
 =\ell_\ominus(w)-(\nu-1)\Delta(w)\quad(\nu=1,2,3,\ldots),\end{array}$$
and this requirers $\Delta(w)\le 0$ for {\it any}
$\alpha\notin\frac 12+\Z$ and {\it any} solution of the first
kind, hence $\Delta(w)\le -1$ by Proposition \ref{LPLUSMINUS}.
Replacing $w$ by $-w$ and $\alpha$ by $-\alpha$, however, we
obtain $\Delta(w)\ge 1$; this contradiction proves the first part.

\medskip To deal with the case $\alpha\in\frac 12+\Z$ it suffices to
consider $\alpha=1/2$. If $w$ is a non-Airy solution of the first
kind the above method applies ``to the right''; we obtain in the
same manner $\Delta(w)\le -1$. Instead of working with $-w$ (which
is prohibited) we now apply the special B\"acklund transformation
$w(z)=-\frac d{dz}\log y(-2^{-1/3}z)$ and obtain from Nevanlinna's
First Main Theorem
 $$\begin{array}{rcl}N_\ominus(r,w)-N_\oplus(r,w)&=&N(2^{-1/3}r,1/y)-N(2^{-1/3}r,y)\cr
 &=&m(2^{-1/3}r,y)-m(2^{-1/3}r,1/y)+O(1)\le O(\log r),\end{array}$$
hence $\Delta(w)\ge 0$. We have thus proved that all sub-normal
solutions of the first kind are Airy solutions (and exist only if
$\alpha\in\frac 12+\Z$).

\section{\bf  Proof of Theorem \ref{HEUREKA}: Second Kind Solutions}\label{HEUREKA2}

Let $w$ be any sub-normal solution of the second kind to equation
[II$_\alpha$]. For $2\alpha\in\Z$ it suffices to consider the
cases $\alpha=0$ and $\alpha=1/2$:

 \begin{itemize}\item If
[II$_{0}$] had a sub-normal solution $w$ of the second kind, then
 $$w_1(z)=-\frac d{dz}\log w(-2^{-1/3}z)$$
were a non-Airy solution of the first kind to equation [II$_{\frac
12}$]. \item Conversely, if [II$_{\frac 12}$] had a sub-normal
solution $w$ of the second kind, then
 $$w_1^2(z)=-2^{1/3}(w'(-2^{1/3}z)-z/2-w^2(-2^{1/3}z)),$$
were a solution of the first kind to equation [II$_{0}$].
\end{itemize}

We now assume $2\alpha\ne\Z$, and set again
 \be{DEV}V=w'+w^2+z/2,\quad w_1=-w-\ds\frac{\alpha+1/2}V.\ee
Then $w_1$ solves [II$_{\alpha+1}$], and the poles of $w$ and
$w_1$ and the zeros of $V$ are related as in (\ref{POLE}).
Associated with any string of poles $(q_k)$ of $w$ is a string of
poles $(p_k)$ of $w_1$ as follows: $p_k=q_k$ if $\res\limits_{q_k}
w=-1$, while $q_{k+1}$ is replaced by $p_{k+1}$ with
$V(p_{k+1})=0$ and $\res\limits_{p_{k+1}}w=0,$ if
$\res\limits_{q_{k+1}} w=-1$; we may assume that this happens for
$k$ even. Since the string $(p_k)$ is already determined by the
sub-string $(p_{2k})$, $p_{2k+1}$ is very close to $q_{2k+1}$, as
is displayed below; $\oplus\!\!\!\circledcirc$ denotes an ``almost
double'' zero of $V$: $V(p_{2k+1})=V(q_{2k+1})=0$ with
$\res\limits_{q_{2k+1}}w=1$ and $\res\limits_{p_{2k+1}}w=0$, while
$\ominus$ and $\oplus$ denote poles with residues $-1$ and $1$,
respectively, for both $w$ (left) and $w_1$ (right); note the
difference to (\ref{VERTEILUNG}).

$$\begin{array}{ccc}\oplus\!\!\!\!\circledcirc\quad\ominus\quad
\oplus\!\!\!\!\circledcirc\quad\ominus\quad\oplus\!\!\!\!\circledcirc\quad\ominus&\qquad&
\ominus\quad\oplus\quad\ominus\quad\oplus\quad\ominus\quad\oplus\end{array}$$

We write $p=p_{2k+1}$ and $q=q_{2k+1}$ and insert
 $$w(z)=b+w'(p)(z-p)+\ts\frac12w''(p)(z-p)^2+\cdots\quad(b=w(p))$$
with $w'(p)=-b-p/2$ and $w''(p)=\alpha+bp+2b^3$ into the
definition (\ref{DEV}) of $w_1$ to obtain (``with a little help
from my friends''---computer algebra software)
 $$w_1(z)=-\frac 1{z-p}-
 \frac{8bp^2+(40b^3-3)p+(48b^5-4b^2+12\alpha b^2)}{6(2\alpha+1)}(z-p)+\cdots$$
Comparing with $w_1(z)=\ds-\frac 1{z-p}+\ds\frac p6(z-p)+\cdots$
yields
 $$p=\frac{1-\alpha-20b^3\pm\sqrt{(1-\alpha)^2-8(1+7\alpha)
 b^3+16b^6}}{8b}\asymp |b|^2\quad(b\to\infty).$$
The special solution $\w=\pm i/\sin(\i\z)$ in (\ref{SPECSOL})
satisfies $|\w(\z)|\ge 2\kappa|\z|^{-1}$ on $|\z|<\delta$ (for
some $\kappa\ge 1$), hence $|w(z)|\ge \kappa|z-q|^{-1}$ on
$|z-q|<\delta|q|^{-1/2}$ if $|q|\ge r_0$. Since
$p-q=o(|p|^{-1/2})$ we obtain $|p|^{-1/2}w(p)\to\infty$ as
$p\to\infty$; this, however, contradicts $|p|\asymp|b|^2$, and
Theorem \ref{HEUREKA} is completely proved. \Ende

\section{\bf Proof of Corollary \ref{deficiency0}}\label{DEFEKT0}

\medskip We have just to consider
solutions of order $\varrho=3$. From the special B\"acklund
transformation $w_1(z)=-\frac d{dz}\log w(-2^{-1/3}z)$, the
estimate
 $$N_\oplus(r,w_1)-N_\ominus(r,w_1)=O(r^{3/2})$$
(see \cite{NSt4}, Thm.\ 6.2), and Nevanlinna's First Main Theorem
we obtain
 $$\begin{array}{rcl}m(r,1/w)&=&T(r,w)-N(r,1/w)+O(1)\cr
 &=&N_\oplus(2^{-1/3}r,w_1)+O(\log r)-N_\ominus(2^{-1/3}r,w_1)=O(r^{3/2}).\end{array}$$
Hence zero is non-deficient for $w$. \Ende

\section{\bf Painlev\'e's Fourth Transcendents: An Outlook}\label{OUTLOOK}

The solutions to Painlev\'e's fourth differential equation
 $$2ww''={w'}^2+3w^4+8zw^3+4(z^2-\alpha)w^2+2\beta\leqno{\rm
 [IV_{\alpha,\beta}]}$$\def\PIV{$[{\rm IV}_{\alpha,\beta}]$}%
are either rational or transcendental meromorphic functions of
order $\varrho$, $2\le\varrho\le 4$. From
 $${w'}^2=w^4+4zw^3+4(z^2-\alpha)w^2-2\beta-4wW\quad(W'=w^2+2zw),$$
follows
 $$w''=2w^3+6zw^2+4(z^2-\alpha)w-2W.$$

\medskip{\bf Rescaling.} The family $(w_h)_{|h|\ge 1}$ of functions
$w_h(\z)=h^{-1}w(h+h^{-1}\z)$ is normal, and every limit function
solves
 $$2\w\w''=\w'^2+3\w^4+8\w^3+4\w^2,$$
with constant solutions $\w=0, -2/3, -2$, and also
 \be{IVre-scaled1}\w'^2=\w^4+4\w^3+4\w^2+4c\w\ee
and
 $$\w''=2\w^3+6\w^2+4\w+2c$$
with $-c$ in the cluster set $\CL_\varepsilon$ of $z^{-3}W(z)$,
which consists of all limits
 $$\lim\limits_{h_n\to\infty}h_n^{-3}W(h_n)\quad(\inf\limits_{n}|h_n|\dist(h_n,\P)\ge\varepsilon).$$
Like in case $\PII$, $\P$ denotes the sequence of non-zero poles
of $w$, and like there it turns out that $\CL=\CL_\varepsilon$ is
independent of $\varepsilon$.

\medskip{\bf Laurent series expansion about poles.} Similar to the
case $\PII$ we have
 $$\begin{array}{rcl}
 w(z)&=&\pm(z-p)^{-1}-p\pm\ts \frac
 13(p^2+2\alpha\mp4)(z-p)+\mathbf{h}(z-p)^2+\cdots\cr
 W(z)&=&-(z-p)^{-1}+[2\mathbf{h}+2(\alpha\mp1)p]+\ts\frac13(4\alpha-p^2\mp
 2)(z-p)+\cdots,\end{array}$$
and the limits
$-c=\lim\limits_{p_n\to\infty}p_n^{-3}[2\mathbf{h}(p_n)+2(\alpha\mp
1)p_n]=\lim\limits_{p_n\to\infty}2p_n^{-3}\mathbf{h}(p_n)$ belong
to $\CL$.

\medskip{\bf Weber-Hermite Solutions.} The r\^{o}le of the Riccati
equations (\ref{Airy}) is taken by the so-called {\it
Weber-Hermite equations}
 \be{HermiteWeber}w'=-2\pm(w^2+2zw-2\alpha).\ee
Their solutions have order of growth $\varrho=2$ and solve
equation [IV$_{\alpha,-2(1\pm\alpha)^2}$]. However, the situation
is more complicated than in case [II$_\alpha$], since their are
several continuous one-parameter families of solutions that can be
reduced to the Weber-Hermite equation, see \S 25 in \cite{GLS}.
They occur for parameters $\beta=-2(2n-1\pm\alpha)^2$ and
$\beta=-2n^2$, respectively; $\alpha$ is arbitrary, and in both
cases $n$ is any integer.

\medskip{\bf Sub-normal solutions.} Our focus
is on the fourth Painlev\'e transcendents with counting function
of poles $n(r,w)=O(r^2).$ The right hand side of
(\ref{IVre-scaled1}) has discriminant $c^3(27c-8)$. It is quite
plausible to analyse the following cases:

\begin{itemize}\item[]{\bf First kind.} $\CL=\{0\}$, $\w'^2=\w^2(\w+2)^2$ with solutions
 $$\w=\ds\frac{2e^{\pm
 2\z+\tau}}{1-e^{\pm 2\z+\tau}}\quad{\rm and}\quad\w=0, -2.$$
The strings of poles $(p_k)_{k=0,1,\ldots}$ are defined by
 $$\begin{array}{rcl}p_{k+1}&=&p_k+\pi i p_k^{-1}+o(|p_k|^{-1})
 \cr p_{k+1}&=&p_k-\pi i p_k^{-1}+o(|p_k|^{-1}),\end{array}$$
hence $p_k\sim \pm(1\pm i)(\pi\,k)^{1/2}$, with counting function
$n(r)\sim \ds\frac{r^2}{2\pi}$;
 $\res\limits_{p_k}w$ is constant on every string. The Weber-Hermite solutions are of the first kind.
 \item[]
 \item[]{\bf Second kind.} $\CL=\{-8/27\}$, $\w'^2=\frac 1{27}\w(3\w+8)(3\w+2)^2$ with solutions
 $$\w=\ds\frac{8}{9\tan^2(\z/\sqrt{3}+\tau)-3}$$
(substitute $3+8/\w=\y^2$) and $\w=-2/3$ (and neither $\w=0$ nor
$\w=-8/3$ occur as limit functions). The strings of poles are
defined by
 $$\begin{array}{rcl}p_{k+1}&=&p_k+\sqrt{3}\pi p_k^{-1}+o(|p_k|^{-1})\cr
 p_{k+1}&=&p_k-\sqrt{3}\pi p_k^{-1}+o(|p_k|^{-1}),\end{array}$$
hence $p_k\sim\pm(2\sqrt{3}\pi\,k)^{1/2}$ and $p_k\sim\pm
i(2\sqrt{3}\pi\,k)^{1/2}$, respectively, with counting function
$n(r)\sim \ds\frac{r^2}{2\sqrt{3}\pi}$; the residues alternate,
$\res\limits_{p_{k+1}}w=-\res\limits_{p_k}w$.
\end{itemize}

\medskip{\bf Yosida Solutions.} For $0,
-8/27\not\in\CL$ all limit functions are elliptic, hence $w$
belongs to the Yosida Class $\Y_{1,1}$ and, in particular,
satisfies $T(r,w)\asymp r^4$. The latter remains true if we
restrict the cluster set of $z^{-3}W(z)$ to any sector.

\bigskip
\small\it Norbert Steinmetz\\ Institut f\"ur Mathematik\\
Technische Universit\"at Dortmund\\ D-44221 Dortmund, Germany\\
E-mail: stein@math.tu-dortmund.de\\
Web: http://www.mathematik.tu-dortmund.de/steinmetz/

\begin{thebibliography}{99}
\bibitem{PB1} P.~Boutroux, Sur
quelques propri\'et\'es des fonctions enti\`eres, {\it Acta.\ Math.}
{\bf 28}, 97-224 (1904).
\bibitem{PB2} P.~Boutroux, Recherches sur
les transcendentes de M. Painlev\'e et l'\'etude asymptotique des
\'equations diff\'erentielles du seconde ordre, {\it Ann.\ \'Ecole
Norm.\ Sup\'er.}\ {\bf 30}, 255-375 (1913) and {\it Ann.\ \'Ecole Norm.\
Sup\'er.}\ {\bf 31}, 99-159 (1914).
\bibitem{GLS} V.~Gromak, I.\ Laine, and S.~Shimomura, {\it Painlev\'e differential equations in the complex plane},
W.\ de Gruyter 2002.
\bibitem{GS} G.\ Gundersen and E.\ Steinbart, A generalization of
the Airy integral for $f''+z^nf=0$, {\it Trans.\ Amer.\ Math.\
Soc.} {\bf 337} (1993), 737-755.
\bibitem{WH} W.K.~Hayman, {\it Meromorphic functions}, Clarendon Press, Oxford 1964.
\bibitem{EH} E.\ Hille, Ordinary differential equations in the
complex domain, Dover Publ. 1997.
\bibitem{HL1} A.\ Hinkkanen and I.\ Laine, Growth results for Painlev\'e transcendents,
{\it Math.\ Proc.\ Cambridge Philos.\ Soc.} {\bf 137} (2004), 645-655.
\bibitem{HL2} A.\ Hinkkanen and I.\ Laine, Growth of second Painlev\'e transcendents, preprint, 46 p. (2011).
\bibitem{RN} R.~Nevanlinna, {\it Eindeutige analytische Funktionen}, Springer 1936.
\bibitem{SSh1} S.\ Shimomura, Growth of the first, the second and the fourth Painlev\'e transcendents,
{\it Math.\ Proc.\ Cambridge Philos.\ Soc.} {\bf 134} (2003), 259-269.
\bibitem{SSh2} S.\ Shimomura, Lower estimates for the growth of the fourth and the second Painlev\'e transcendents,
{\it Proc.\ Edinb.\ Math.\ Soc.} {\bf 47} (2004), 231-249.
\bibitem{NSt1} N.\ Steinmetz, On Painlev\'e's equations I, II and IV, {\it J.\
d'Analyse Math.} {\bf 82} (2000), 363-377.
\bibitem{NSt2} N.~Steinmetz, Value distribution of the Painlev\'e transcendents,
{\it Israel J.\ Math.} {\bf 128} (2002), 29-52.
\bibitem{NSt4} N.\ Steinmetz, Global properties of the Painlev\'e transcendents. New results and open questions,
{\it Ann.\ Acad.\ Sci.\ Fenn.\ A I Math.} {\bf 30} (2005), 71-98.
\bibitem{NSt3} N.~Steinmetz, The Yosida class is universal, preprint, 14 p. (2011).
\end{thebibliography}
\end{document}